\documentclass[11pt, reqno]{amsart}

\usepackage{amssymb,latexsym,amsmath,amsfonts}
\usepackage{enumitem}
\usepackage{marvosym}
\usepackage{mathrsfs}
\usepackage{graphicx}
\usepackage{hyperref}
\usepackage[usenames]{color}
\usepackage{orcidlink}
\usepackage{tikz-cd}

\numberwithin{equation}{section}

\theoremstyle{definition}
\newtheorem{definition}{Definition}[section]

\theoremstyle{remark}
\newtheorem{remark}[definition]{Remark}
\theoremstyle{plain}
\newtheorem{theorem}[definition]{Theorem}
\newtheorem{result}[definition]{Result}
\newtheorem{lemma}[definition]{Lemma}
\newtheorem{proposition}[definition]{Proposition}

\definecolor{DPurple}{rgb}{0.46,0.2,0.69}

\newcommand{\supp}{{\rm supp}}

\usepackage[backend=biber,style=alphabetic]{biblatex}
\begin{document}

\title{ A note on the product of holomorphic correspondences}
\author{Bharath Krishna Seshadri \orcidlink{0009-0004-0541-8040} \ \and \ Sathi Trikkadeeri \orcidlink{0009-0000-1444-4604}}

\address{Indian Institute of Science Education and Research Thiruvananthapuram}
\email{bharathmaths21@iisertvm.ac.in}

\address{Indian Institute of Science Education and Research Thiruvananthapuram}
\email{sathitm23@iisertvm.ac.in}

\begin{abstract}
    In this article we study the concept of product of two holomorphic correspondences which was introduced in \cite{sabha26}. Beginning with some basic results regarding the orbit space associated with the product correspondence, we then move on to a detailed analysis on the thermodynamic formalism of the same. Some of the results look at the behaviour of the product correspondence with respect to the individual correspondences. Finally, we explore different families of correspondences with expanding nature.   
\end{abstract}

\keywords{Product correspondence, Dinh Sibony measure, Entropy, Pressure, Unique equilibrium state, Distance expanding correspondence}
\subjclass[2020]{37F05, 37A35, 37F44}

\maketitle

\section{Introduction}
\label{sec:one}

This work is intended to provide an extensive study on the concept of product of holomorphic correspondences, which was discussed in \cite{sabha26}. This notion was considered in \cite{sabha26} with regards to proving some results on weak mixing. Our main aim in this work is to analyse the dynamics of the product correspondence with respect to the dynamics of the individual correspondences that make up the product. This is again inspired by the scenario in the study of single valued maps where by considering the product of two maps, the connection between the individual and combined behaviour of two different dynamical systems are studied. Many results on this topic are well known in the literature for the study of maps.  \vspace{2mm}

The article shall mainly focus on studying the ergodic theory of product of two holomorphic correspondences, each of which is defined on a compact K\"{a}hler manifold. A special focus is given to the theory of thermodynamic formalism. We also build upon recent works in this field and study different classes of correspondences, namely forward expansive, expansive, distance expanding, globally distance expanding, and restrictive holomorphic correspondences. A key question in this regard is about which properties are preserved under the product operation. We now introduce some basic definitions involved in our work.

\begin{definition}
Let $X$ be a compact K\"{a}hler manifold of dimension $k$. Also, let $\pi_1$ and $\pi_2$ denote the projections from $X \times X$ onto the first and second coordinates, respectively. Consider a finite linear combination of irreducible complex sub varieties of $X \times X$ (with each of them having dimension $k$), given by 
\begin{equation} 
\label{eqn:holocorr} 
\Gamma  = \sideset{}{'} \sum \limits_{j\, =\, 1}^{N} \Gamma_{j}, 
\end{equation}
where $\sideset{}{'} \sum$ accounts for the fact that sub varieties are allowed to repeat in the sum. $\Gamma$ is said to be a \emph{holomorphic correspondence} on $X$, when $\pi_i|_{\Gamma_j}$ is surjective, finite-to-one map, for each $i=1,2$ and for all $1 \leq j \leq N$.
\end{definition}

The \emph{adjoint} of a holomorphic correspondence $\Gamma$ is given by $\Gamma^\dagger = \sideset{}{'} \sum \limits_{j=1}^N \Gamma_j^\dagger$, where $\Gamma_j^\dagger = \{ (x, y) \in X \times X : (y, x) \in \Gamma_{j} \}$. On any subset $A \subseteq X$, the set-valued maps induced from $\Gamma$ and $\Gamma^{\dagger}$, to be denoted by $F$ and $F^\dagger$, act as $F(A) := \pi_{2} (\pi_{1}^{-1} (A) \cap |\Gamma|)$ and $F^{\dagger}(A) := \pi_{1} (\pi_{2}^{-1} (A) \cap |\Gamma|)$, where $|\Gamma| := \bigcup \limits_{j=1}^N \Gamma_j$. These yield the `image set' $F(A)$ and the `preimage set' $F^\dagger(A)$ of $A$. For generic $x \in X$, the cardinalities (counting in multiplicity) of $F(x)$ and $F^\dagger(x)$ are called the \emph{forward degree} (denoted by $d_f$) and \emph{backward degree} (denoted by $d_t$) respectively. We will be working with Borel probability measures throughout this article. We say a Borel probability measure $\mu$ on $X$ is \emph{$F^*$-invariant} if it satisfies $F^{*}\mu = d_t\mu$, where $F^{*}\mu$ denotes the pullback of the measure $\mu$ by the correspondence $F$. \vspace{2mm}

Now we define the product of two holomorphic correspondences, which has been developed and studied in detail in  \cite{sabha26}. Let $X_1$ and $X_2$ be two compact K\"{a}hler manifolds of dimensions $k_1$ and $k_2$ respectively and $\Gamma_1  =\displaystyle \sideset{}{'} \sum\limits_{1\leq i\leq N_1} \Gamma_i^1$ and $\Gamma_2 = \displaystyle \sideset{}{'} \sum\limits_{1\leq j\leq N_2} \Gamma_j^2$ be two holomorphic correspondences defined on $X_1$ and $X_2$. Then the \emph{product correspondence} is defined as  
\begin{equation}
\label{eqn:proddefn}
     \Gamma_{1} \times \Gamma_{2} := \sideset{}{'} \sum \limits_{1\leq i\leq N_1} \, \sideset{}{'} \sum \limits_{1\leq j\leq N_2} S^{-1} (\Gamma_i^{1} \times \Gamma_j^{2}),
\end{equation}
where $S: X_1 \times X_2 \times X_1 \times X_2 \longrightarrow X_1 \times X_1 \times X_2 \times X_2$ is the holomorphic function given by $S(x_1, x_2, x_3, x_4) = (x_1, x_3, x_2, x_4)$. It is easy to verify that $\Gamma_{1} \times \Gamma_{2}$ as defined in Equation \eqref{eqn:proddefn} is a holomorphic correspondence on $X_1 \times X_2$. The set valued map corresponding to $\Gamma_1 \times \Gamma_2$, denoted by $F_1 \times F_2$ is given by $(F_1 \times F_2) (x_1, x_2) = F_1 (x_1) \times F_2 (x_2)$. Now, if $d_{t_1}$ and $d_{t_2}$ are the topological degrees of $F_1$ and $F_2$, respectively, then the topological degree of $F_1\times F_2 = d_{t_1}d_{t_2}$. Also if $\mu_1$ is $F_{1}^{*}$ invariant and $\mu_2$ is $F_2^{*}$ invariant, then the product measure $\mu_1 \times \mu_2$ is $(F_1 \times F_2)^{*}$ invariant (see Proposition 5.1 in \cite{sabha26}). All the ideas discussed above can be found in detail in the literature. In particular, one may refer to \cite{bs_p1}, \cite{londhe_ergodicity}, and \cite{sabha26}.\vspace{2mm}

The paper is organised as follows. In Section \ref{sec:prelims} we define and study the orbit space associated with a holomorphic correspondence while in Section \ref{sec:3} we prove some foundational results connecting the orbit spaces associated with the product correspondence and that of the individual holomorphic correspondences. Section \ref{sec:4} is devoted for the study of thermodynamic formalism of the product correspondences. In the case of maps we have entropy of the product of two maps to be equal to the sum of their individual entropies. But Theorem \ref{thm:entropycorrineq} shows that it is an inequality in our situation instead. As for topological pressure, Theorem \ref{thm:presseqcorr} proves that the analogous equality in case of maps still holds in our system, that is ${\rm Pr} (\Gamma_1 \times \Gamma_2, f_1 \times f_2)={\rm Pr} (\Gamma_1, f_1)+{\rm Pr} (\Gamma_2, f_2)$. The last section introduces distance expanding and forward expansive correspondences and shows that their respective products also inherit their behaviours (see Proposition \ref{prop:proddistexp}). Another important result we have in this section is regarding the unique equilibrium state. In Theorem \ref{thm:uniqequi} we prove that, under the dynamics of product of two distance expanding holomorphic correspondences $F_1$ and $F_2$, the unique equilibrium state associated to the potential $f_{1}\times f_{2}$ is given by the product of the unique equilibrium states associated to the individual potentials. Results \ref{thm:globexp} and \ref{rslt:GDE-Exp} discuss the relationships between the different classes of holomorphic correspondences that we are interested in.  

\section{The space of orbits}
\label{sec:prelims}
In Section \ref{sec:one} we saw how the dynamics of a holomorphic correspondence works when it is considered as a set valued map. We now construct an alternate framework that will provide a convenient perspective in some contexts for examining the dynamics of the holomorphic correspondence $\Gamma$ on $X$. The idea of permissible paths under the iteration of the correspondence has been considered in \cite{DinhSibony2008}, \cite{bs_entropy} and \cite{Sub2026}, to quote a few works from the literature. We will mostly follow the methods given in \cite{Sub2026}. \vspace{2mm}

Given $n \in \mathbb{Z}_{+}$, we denote the set of all permissible forward orbits of iteration, of length $n$ by 
\begin{eqnarray*} 
\mathscr{P}_{n}^{\Gamma} \left( X \right) & = & \Bigg\{ \left( x_{0}, x^{(1)}_{j_{1}}, \cdots, x^{(n)}_{j_{n}};\; \alpha_{1}, \cdots, \alpha_{n} \right)\ :\ \left( x^{(r - 1)}_{j_{r - 1}}, x^{(r)}_{j_{r}} \right) \in \Gamma_{\alpha_{r}}\ \text{where} \ x^{(0)}_{j_{0}} = x_{0} \in X, \\
& &  \alpha_{r} \in \left\{ 1, 2, \cdots, N \right\} \text{and}\ 1 \le j_{r} \le \lambda_{\alpha_{r}} \left( x^{(r - 1)}_{j_{r - 1}} \right), \ \text{for}\ 1 \le r \le n \Bigg\}, 
\end{eqnarray*} 
where $\lambda_{\alpha_{r}} \left( x^{(r - 1)}_{j_{r - 1}} \right) = \# \left( F \left (x_{j_{r - 1}}^{(r - 1)} \right) \cap \pi_{2} (\Gamma_{\alpha_{r}})\right)$. Similarly, we define the space of all possible infinitely long orbits starting at an arbitrary point in $X$ as follows.

\begin{eqnarray*} 
\mathscr{P}^{\Gamma} \left( X \right) & = & \Bigg\{ \left( x_{0}, x^{(1)}_{j_{1}}, x^{(2)}_{j_{2}}, \cdots;\; \alpha_{1}, \alpha_{2}, \cdots \right)\ :\ \left( x^{(r - 1)}_{j_{r - 1}}, x^{(r)}_{j_{r}} \right) \in \Gamma_{\alpha_{r}}\ \text{where} \ x^{(0)}_{j_{0}} = x_{0} \in X, \\ & & \alpha_{r} \in \left\{ 1, 2, \cdots, N \right\} \text{and}\ 1 \le j_{r} \le \lambda_{\alpha_{r}} \left( x^{(r - 1)}_{j_{r - 1}} \right), \ \text{for}\ r \in \mathbb{Z}_{+} \Bigg\}. 
\end{eqnarray*} 

We denote a point in $\mathscr{P}^{\Gamma} \left( X \right)$ by $\mathfrak{X}^{+} \left( x_{0}; \boldsymbol{\alpha} \right)_{\boldsymbol{j}}$, where $x_{0} \in X$ is the starting point of the orbit,  $\boldsymbol{\alpha} = \left( \alpha_{1}, \alpha_{2}, \cdots \right)$ tracks the variety through which the orbit is moving forward and $\boldsymbol{j} = \left( j_{1}, j_{2}, \cdots \right)$ distinguishes the images of the same point that are found in the same sub variety. Similarly, a point in $\mathscr{P}_{n}^{\Gamma} \left( X \right)$ will be denoted by $\mathfrak{X}_{n}^{+} (x_{0}; \boldsymbol{\alpha})_{\boldsymbol{j}}$. For $0 \le r \le n$ we define, $\Pi_{(r,n)}^{\Gamma} : \mathscr{P}_{n}^{\Gamma} (X) \to X \ \ \text{to be} \ \ \Pi_{(r,n)}^{\Gamma} (\mathfrak{X}_{n}^{+} (x_{0}; \boldsymbol{\alpha})_{\boldsymbol{j}})  =   x_{j_r}^{(r)}$ and for $1 \le r \le n$ define, $\text {Proj}_{(r,n)}^{\Gamma}: \mathscr{P}_{n}^{\Gamma} (X) \to \left\{ 1, 2, \cdots, N \right\} \ \ \text{to be} \ \ \text {Proj}_{(r,n)}^{\Gamma} (\mathfrak{X}_{n}^{+} (x_{0}; \boldsymbol{\alpha})_{\boldsymbol{j}})  = \alpha_{r}.$ Note that we can similarly define projections $\Pi_{r}^{\Gamma}$ for any $r \ge 0$ and $\text{Proj}_{r}^{\Gamma}$ for all $r \ge 1$ on the space $\mathscr{P}^{\Gamma} \left( X \right)$ as well with appropriate modifications. \vspace{2mm}

The shift map $\sigma^{\Gamma}$ associated with the holomorphic correspondence, defined on $\mathscr{P}^{\Gamma} \left( X \right)$ is given by $\sigma^{\Gamma} \left( \mathfrak{X}^{+} \left( x_{0}; \boldsymbol{\alpha} \right)_{\boldsymbol{j}} \right) = \left( x^{(1)}_{j_{1}}, x^{(2)}_{j_{2}}, \cdots;\; \alpha_{2}, \alpha_{3},  \cdots \right) = \mathfrak{X}^{+} \left( x^{(1)}_{j_{1}}; \sigma \boldsymbol{\alpha} \right)_{\sigma \boldsymbol{j}}$, where $\sigma \boldsymbol{\alpha}$ and $\sigma \boldsymbol{j}$ respectively denote the infinite-lettered word got by removing the first letter of $\boldsymbol{\alpha}$ and $\boldsymbol{j}$. We also define the metric $d$ on $\mathscr{P}^{\Gamma} \left( X \right)$ given by,
\[
d \left( \mathfrak{X}^{+} \left( x_{0}; \boldsymbol{\alpha} \right)_{\boldsymbol{j}}, \mathfrak{X}^{+} \left( y_{0}; \boldsymbol{\beta} \right)_{\boldsymbol{k}} \right) = \max \left\{ \sum \limits_{n \ge 0} \frac{d_{X}\left( x_{j_{n}}^{(n)},  y_{k_{n}}^{(n)}\right)}{2^{n} \left(1+ d_{X}\left( x_{j_{n}}^{(n)},  y_{k_{n}}^{(n)}\right)\right)},\;\frac{1}{2^{\min\{n \mid \alpha_n \ne \beta_n \}}} \right\}, 
\]
where $d_{X}$ denotes the distance function on the K\"{a}hler manifold $X$. Observe that $\mathscr{P}^{\Gamma} \left( X \right)$ is a compact metric space and the map $\sigma^{\Gamma}$ is continuous with respect to this metric. The metric defined above induces the product topology on $\mathscr{P}^{\Gamma} \left( X \right)$ that is obtained as a subspace of $X^{\mathbb{N}} \times \{ 1, 2, \cdots, N \} ^ \mathbb{N}$.\vspace{2mm}

Let $\Gamma_1$ and $\Gamma_2$ be holomorphic correspondences on the compact K\"{a}hler manifolds $X_1$ and $X_2$ respectively and let $\Gamma_{1} \times \Gamma_{2}$ denote their product. An element of $\mathscr{P}^{\Gamma_{1} \times \Gamma_{2}} \left( X_1 \times X_2 \right)$ is given by $\left( (x_{0},y_{0}), (x^{(1)}_{j_{1}},y^{(1)}_{k_{1}}), \cdots;\; (\alpha_{1},\beta_{1}), (\alpha_{2}, \beta_{2}) \cdots  \right)$ and is denoted by $\mathfrak{X}^{+} \left((x_{0}, y_{0});(\boldsymbol{\alpha} \times \boldsymbol{\beta})\right)_{(\boldsymbol{j} \times{\boldsymbol{k}})}$. Here the tuple $(\alpha_j, \beta_k)$ represents the variety $S^{-1} (\Gamma_{\alpha_j} \times \Gamma_{\beta_{k}})$ as mentioned in Equation \eqref{eqn:proddefn}. Moreover, for $\boldsymbol{\alpha} = (\alpha_1, \alpha_2, \cdots) \ \text{and} \ \boldsymbol{\beta} = (\beta_1, \beta_2, \cdots)$, the symbol $(\boldsymbol{\alpha} \times \boldsymbol{\beta})$ denotes infinite tuple given by $\left( (\alpha_1, \beta_1), (\alpha_2, \beta_2), \cdots \right)$. The symbol $(\boldsymbol{j} \times \boldsymbol{k})$ must be interpreted similarly. \vspace{2mm}

We end this section by introducing some measure theoretic ideas on the orbit space associated to a correspondence. Let $\mathscr{M} \left( \mathscr{P}^{\Gamma} \left( X \right) \right)$ denote the collection of all Borel probability measures supported on $\mathscr{P}^{\Gamma} \left( X \right)$. The space $\mathscr{M} \left( \mathscr{P}^{\Gamma} \left( X \right) \right)$ is compact under the weak* topology, and this topology is metrizable. Let $\mathscr{M} \left( \mathscr{P}^{\Gamma} \left( X \right), \sigma^{\Gamma} \right)$ represent the subset of $\mathscr{M} \left( \mathscr{P}^{\Gamma} \left( X \right) \right)$ consisting of all probability measures that are invariant under $\sigma^{\Gamma}$. The space $\mathscr{M} \left( \mathscr{P}^{\Gamma}(X), \sigma^{\Gamma} \right)$ is non-empty and compact under the weak* topology. Further, the projection $\Pi_0^{\Gamma}: \mathscr{P}^{\Gamma} \left( X \right) \to X$ given by $\Pi_{0}^{\Gamma} (\mathfrak{X}^{+} \left( x_{0}; \boldsymbol{\alpha} \right)_{\boldsymbol{j}}) = x_0$, induces a push forward map, $(\Pi_0^{\Gamma})_{*}$ from $\mathscr{M} \left( \mathscr{P}^{\Gamma} \left( X \right)\right)$ to $\mathscr{M} \left( X \right)$. Define $\mathscr{S}^{\Gamma} \subset \mathscr{M} \left( X \right)$ as, $\mathscr{S}^{\Gamma} := (\Pi_0^{\Gamma})_{*} \left(\mathscr{M} \left( \mathscr{P}^{\Gamma} \left( X \right), \sigma^{\Gamma} \right)\right)$. We shall extensively work with $\mathscr{S}^{\Gamma}$ in the upcoming sections. 

\section{Some results on the dynamics of product correspondence}
\label{sec:3}
After having set up the basic notations and the dynamical systems involved in this paper we now turn our attention towards proving some basic properties related to the dynamics of the product correspondence. We mainly work with the space of orbits in this section and also occasionally go to the ambient K\"{a}hler manifold space. We prove some conjugacy relations that establish an isomorphism between two dynamical systems. This shall be crucial in our quest to understand the product correspondences better and also make the calculations a lot simpler. In what follows, for maps $f_1$ and $f_2$ defined on sets $S_1$ and $S_2$, the notation $(f_1, f_2)$ denotes the map on $S_1 \times S_2$ given by $(f_1, f_2) (s_1, s_2) = (f_1(s_1), f_2(s_2))$.

\begin{lemma}
The orbit spaces $\mathscr{P}^{\Gamma_1}{(X_1)} \times \mathscr{P}^{\Gamma_2}{(X_2)}$ and $\mathscr{P}^{{\Gamma_1}\times {\Gamma_2}}{(X_1 \times X_2)}$ are homeomorphic. Here, the space $\mathscr{P}^{\Gamma_1}{(X_1)} \times \mathscr{P}^{\Gamma_2}{(X_2)}$ is assigned the product topology.
\end{lemma}

\begin{proof}
Consider the point $\mathfrak{X}^+(x_0; \boldsymbol{\alpha})_{\boldsymbol{j}} = \left( x_{j_0}^{(0)}, x_{j_1}^{(1)}, \dots ;\, \alpha_1, \alpha_2, \dots \right) \in \mathscr{P}^{\Gamma_1}{(X_1)}$ and the point $\mathfrak{X}^+(y_0; \boldsymbol{\beta})_{\boldsymbol{k}} = \left( y_{k_0}^{(0)}, y_{k_1}^{(1)}, \dots ;\, \beta_1, \beta_2, \dots \right) \in \mathscr{P}^{\Gamma_2}{(X_2)}$. Define the map $\Phi$ from the space $\mathscr{P}^{\Gamma_1}{(X_1)} \times \mathscr{P}^{\Gamma_2}{(X_2)}$ to $\mathscr{P}^{{\Gamma_1}\times{\Gamma_2}}{(X_1 \times X_2})$ as follows.
\begin{eqnarray}
\label{eqn:homeo}
    \Phi\left( \mathfrak{X}^+(x_0;\boldsymbol{\alpha})_{\boldsymbol{j}} ,\, \mathfrak{X}^+(y_0;\boldsymbol{\beta})_{\boldsymbol{k}} \right) & = & \left( \left(x_{j_0}^{(0)}, y_{k_0}^{(0)}\right),\, \left(x_{j_1}^{(1)}, y_{k_1}^{(1)}\right), \dots ;\, (\alpha_1, \beta_1), (\alpha_2, \beta_2), \dots \right) \nonumber \\
    & = & \left( \left( \left( \Pi_{r}^{\Gamma_1}(\mathfrak{X}^+(x_0; \boldsymbol{\alpha})_{\boldsymbol{j}}), \Pi_{r}^{\Gamma_2}(\mathfrak{X}^+(y_0; \boldsymbol{\beta})_{\boldsymbol{k}}) \right) \right)_{r \ge 0}; \right. \nonumber \\
    &  & \hspace{1mm} \left.  \left( \left( \text{Proj}_{r}^{\Gamma_1}(\mathfrak{X}^+(x_0; \boldsymbol{\alpha})_{\boldsymbol{j}}), \text{Proj}_{r}^{\Gamma_2}(\mathfrak{X}^+(y_0; \boldsymbol{\beta})_{\boldsymbol{k}}) \right) \right)_{r \ge 1} \right)
\end{eqnarray}
Since the maps $\Pi_r^{\Gamma_1}$ and $\Pi_r^{\Gamma_2}$ are continuous for all $r \ge 0$. and the maps $\text{Proj}_r^{\Gamma_1}$ and $\text{Proj}_r^{\Gamma_2}$ are continuous for all $r \ge 1$, we have $\Phi$ to be a continuous map. The inverse of $\Phi$ defined on $\mathscr{P}^{{\Gamma_1}\times{\Gamma_2}}{(X_1 \times X_2)}$ is given by, 
\[
\Phi^{-1}\left( \mathfrak{X}^{+} ((x_{0}, y_{0});(\boldsymbol{\alpha} \times \boldsymbol{\beta}))_{(\boldsymbol{j} \times {\boldsymbol{k}})} \right) = \left( \left(x_{j_0}^{(0)}, x_{j_1}^{(1)}, \dots ; \alpha_1, \alpha_2, \dots\right), \left(y_{k_0}^{(0)}, y_{k_1}^{(1)}, \dots ; \beta_1, \beta_2, \dots\right) \right)
\]
An easy computation shows that $\Phi^{-1}$ is continuous as well. Hence the result holds.
\end{proof}

Thus for all $r \ge 0$, we have the following commutative diagram. 

\begin{equation}
\label{dgm:1}
\begin{tikzcd}[row sep = huge, column sep = huge]
\mathscr{P}^{\Gamma_1}{(X_1)} \times \mathscr{P}^{\Gamma_2}{(X_2)} \arrow[r, "{\left(\Pi_r^{\Gamma_{1}}, \Pi_r^{\Gamma_{2}}\right)}"] \arrow[d, "\Phi"'] & X_1 \times X_2 \arrow[d, "\text{Id}"] \\
\mathscr{P}^{{\Gamma_1}\times{\Gamma_2}}{(X_1 \times X_2)} \arrow[r, "{\Pi_r^{\Gamma_1 \times \Gamma_2}}"] & X_1 \times X_2
\end{tikzcd}
\end{equation}

\begin{remark}
\label{rem:projncomm}
    It is easy to see that a commutative diagram along the lines of Diagram \eqref{dgm:1} exists between the spaces $\mathscr{P}^{\Gamma_1}_{n}(X_1) \times \mathscr{P}^{\Gamma_2}_{n}(X_2)$ and $\mathscr{P}^{{\Gamma_1}\times{\Gamma_2}}_{n}{(X_1 \times X_2)}$ as well with the corresponding homeomorphism being denoted by $\Phi_{n}$ and the projections considered being $\Pi_{(r, n)}^{\Gamma_{1} \times \Gamma_{2}}$ and $\text {Proj}_{(r,n)}^{\Gamma_{1} \times \Gamma_{2}}$.
\end{remark}

Now we show that the two spaces are not just homeomorphic, but when they are considered as dynamical systems with appropriate measures and shift maps, they are indeed isomorphic measure preserving systems. Before stating the results, we introduce a few notations. Fix $i = 1, 2$. We have $\mathcal{B}_i$ to denote the Borel $\sigma$ - algebra of $X_i$ and measures $\mu_i$ are Borel probability measures on $\mathscr{P}^{\Gamma_i}{(X_i)}$ that are invariant under the maps $\sigma^{\Gamma_{i}}$. Also, $\mathcal{B}$ denotes the Borel $\sigma-$ algebra on $\mathscr{P}^{{\Gamma_1}\times{\Gamma_2}}{(X_1 \times X_2)}$.

\begin{lemma}
\label{lem:MPS}
The system $\left( \mathscr{P}^{\Gamma_1}{(X_1)} \times \mathscr{P}^{\Gamma_2}{(X_2)},\, (\sigma^{\Gamma_1}, \sigma^{\Gamma_2}), \, \mathcal{B}_1 \times \mathcal{B}_2,\, \mu_1 \times \mu_2 \right)$ and the system $\left( \mathscr{P}^{{\Gamma_1}\times{\Gamma_2}}{(X_1 \times X_2)},\, \sigma^{\Gamma_1 \times \Gamma_2} , \, \mathcal{B},\, \Phi_*(\mu_1 \times \mu_2) \right)$ are isomorphic measure preserving dynamical systems.
\end{lemma}
\begin{proof}
Firstly it is easy to observe that we have the following commutative diagram. 
\[
\begin{tikzcd}[row sep = huge, column sep = huge]
\mathscr{P}^{\Gamma_1}{(X_1)} \times \mathscr{P}^{\Gamma_2}{(X_2)} \arrow[r, "{\left(\sigma^{\Gamma_1}, \sigma^{\Gamma_2}\right)}"] \arrow[d, "\Phi"'] & \mathscr{P}^{\Gamma_1}{(X_1)} \times \mathscr{P}^{\Gamma_2}{(X_2)} \arrow[d, "\Phi"] \\
\mathscr{P}^{{\Gamma_1}\times{\Gamma_2}}{(X_1 \times X_2)} \arrow[r, "{\sigma^{\Gamma_1 \times \Gamma_2}}"] & \mathscr{P}^{{\Gamma_1}\times{\Gamma_2}}{(X_1 \times X_2)}
\end{tikzcd}
\]
We use this to now prove that $\Phi_*(\mu_1 \times \mu_2) \in \mathscr{M} \left( \mathscr{P}^{{\Gamma_1}\times{\Gamma_2}}{(X_1 \times X_2)},\, \sigma^{\Gamma_1 \times \Gamma_2} \right)$. Let $A \subset \mathscr{P}^{{\Gamma_1}\times{\Gamma_2}}{(X_1 \times X_2})$ be a Borel set. Then we have,
\begin{align*}
\left( \Phi_*(\mu_1 \times \mu_2) \right) \left( (\sigma^{\Gamma_1 \times \Gamma_2})^{-1}(A) \right)
& =  (\mu_1 \times \mu_2) \left( \left( \sigma^{\Gamma_1 \times \Gamma_2} \circ \Phi \right)^{-1}(A) \right) \\
& =  (\mu_1 \times \mu_2) \left( \left( \Phi \circ \left(\sigma^{\Gamma_1}, \sigma^{\Gamma_2}\right) \right)^{-1}(A) \right) \\
& =  (\mu_1 \times \mu_2) \left( \left(\sigma^{\Gamma_1}, \sigma^{\Gamma_2}\right)^{-1} \left( \Phi^{-1}(A) \right) \right) \\
&= (\mu_1 \times \mu_2) \left( \Phi^{-1}(A) \right) \\
&= \Phi_*(\mu_1 \times \mu_2)(A).
\end{align*}
As $\Phi$ and $\Phi^{-1}$ are continuous, they are measurable as well. Further, $\Phi$ preserving the measures of the respective systems comes for free due to the choice $\Phi_*(\mu_1 \times \mu_2)$ for the measure on the second system. This completes the proof of the theorem.
\end{proof}

\begin{remark}
\label{rem:entpreseq}
    An immediate consequence of Lemma \ref{lem:MPS} is that the measure theoretic entropy $h_{\mu_1 \times \mu_2} ((\sigma^{\Gamma_1}, \sigma^{\Gamma_2}))$ of $(\sigma^{\Gamma_1}, \sigma^{\Gamma_2})$ with respect to $\mu_1 \times \mu_2$ is equal to $h_{\Phi_*(\mu_1 \times \mu_2)} (\sigma^{\Gamma_1 \times \Gamma_2})$, the measure theoretic entropy of $\sigma^{\Gamma_1 \times \Gamma_2}$ with respect to $\Phi_*(\mu_1 \times \mu_2)$. Similarly, by using Remark \ref{rem:projncomm} and Lemma \ref{lem:MPS} we have that, ${\rm Pr} (\sigma^{\Gamma_1 \times \Gamma_2}, f ) = {\rm Pr} ((\sigma^{\Gamma_1}, \sigma^{\Gamma_2}), f \circ \Phi)$, for all $f \in \mathcal{C} (\mathscr{P}^{{\Gamma_1}\times{\Gamma_2}}{(X_1 \times X_2}) , \mathbb{R})$. One may refer \cite{walters_book} for details regarding the same.
\end{remark}

Having proved some measure theoretic results on the space of orbits, we now turn our attention to the measure theoretic aspects of the dynamics of the correspondence on the K\"{a}hler manifold.
\begin{proposition} 
\label{prop:SGamma}
    Suppose $\nu_1 \in \mathscr{S}^{\Gamma_1}$ and $\nu_2 \in \mathscr{S}^{\Gamma_2}$, then $\nu_1 \times \nu_2 \in \mathscr{S}^{\Gamma_1 \times \Gamma_2}$.
\end{proposition}
\begin{proof}
For $i=1,2$, choose $\mu_i \in \mathscr{M}(\mathscr{P}^{\Gamma_i}{(X_i)}, \sigma^{\Gamma_i})$ such that $\nu_i = (\Pi_0^{\Gamma_i})_{*}\mu_i$. We show that $(\Pi_0^{\Gamma_1 \times \Gamma_2})_* \left( (\Phi)_* (\mu_1 \times \mu_2) \right) = \nu_1 \times \nu_2.$ In order to do so it is enough to prove that,
\begin{equation}
\label{eqn:1}
\int \limits_{X_1 \times X_2} f \, \mathrm{d} \hspace{-1mm}\left( \left(\Pi_0^{\Gamma_{1} \times \Gamma_{2}}\right)_* ((\Phi)_* (\mu_1 \times \mu_2)) \right) = \int \limits_{X_1 \times X_2} f \, \mathrm{d}(\nu_1 \times \nu_2) 
\end{equation}
for all $f \in \mathcal{C}(X_1 \times X_2, \mathbb{R})$. Let $f_1 \in \mathcal{C} (X_1, \mathbb{R})$ and $f_2 \in \mathcal{C} (X_2, \mathbb{R})$. Define $f_1 \times f_2: X_1 \times X_2 \to \mathbb{R}$ given by $(f_1 \times f_2) (x_1, x_2) = f_1(x_1)f_2(x_2)$. Since functions of the form $f_1 \times f_2$ are dense in $\mathcal{C} (X_1 \times X_2, \mathbb{R})$, it is enough to prove Equation \eqref{eqn:1} holds for such functions. Fixing $f = f_1 \times f_2$ as mentioned above and using the Diagram \eqref{dgm:1}, we observe that,
\begin{eqnarray*} 
\int \limits_{X_1 \times X_2} f \, \mathrm{d} \hspace{-1mm}\left( \left(\Pi_0^{\Gamma_{1} \times \Gamma_{2}}\right)_* ((\Phi)_* (\mu_1 \times \mu_2)) \right) 
& = & \int\limits_{\mathscr{P}^{\Gamma_1}{(X_1)} \times \mathscr{P}^{\Gamma_2}{(X_2)}} (f_1 \times f_2) \circ (\Pi_0^{\Gamma_{1}}, \Pi_0^{\Gamma_{2}}) \, \mathrm{d}(\mu_1 \times \mu_2)\\
& = & \int\limits_{\mathscr{P}^{\Gamma_1}{(X_1)}} (f_1 \circ \Pi_0^{\Gamma_1}) \, \mathrm{d}\mu_1 \int\limits_{\mathscr{P}^{\Gamma_2}{(X_2)}} (f_2 \circ \Pi_0^{\Gamma_2}) \, \mathrm{d}\mu_2 \\
& = & \int\limits_{X_1} f_1 \, \mathrm{d} \hspace{-1mm} \left( (\Pi_0^{\Gamma_{1}})_* \mu_1 \right)  \int\limits_{X_2} f_2 \, \mathrm{d} \hspace{-1mm} \left( (\Pi_0^{\Gamma_{2}})_* \mu_2 \right) \\
& = & \int\limits_{X_1 \times X_2} f \, \mathrm{d}(\nu_1 \times \nu_2).
\end{eqnarray*}
\end{proof}
We end this section with a result on the connection between the Dinh Sibony measure of product of two holomorphic correspondences with respect to the individual ones. This measure is very crucial in the study of holomorphic correspondences and it plays a key role in understanding the thermodynamic formalism in this dynamical system. Further, the support of this measure is considered to be an analogue of the Julia set from the case of rational maps. \vspace{2mm}

Dinh and Sibony in \cite{dinhsibony_og} proved that for a holomorphic correspondence $F$  defined on a compact K\"{a}hler manifold $X$ of dimension $k$ satisfying $d_t > d_{k - 1}$, where $d_{k - 1}$ is the dynamical degree of $F$ of order $k - 1$ and $d_t$ is the topological degree of $F$, there exists a pluripolar set $\mathcal{E} \subsetneq X$ such that $\displaystyle{\frac{1}{d^{n}} \left( F^{n} \right)^{*} \delta_{x} \to \mu_{F}}, \text{in the weak* topology,\ for all}\ x \in X \setminus \mathcal{E}.$ Moreover, the measure $\mu_F$ (which we refer to as the Dinh Sibony measure) is $F^{*}$ invariant and $\mu_F(\mathcal{E}) = 0$. Here, $\delta_{x}$ denotes the Dirac delta measure at the point $x$. It must be noted that, an analogous result has been proved for some other families of holomorphic correspondences too under various hypotheses. Some of the works in connection with that include \cite{bs_p1}, \cite{corr_RS}, \cite{hemmingsson} and \cite{de_parra}, thereby resulting in a possibility of studying Dinh Sibony measure for such families too. 

\begin{proposition}
\label{prop:dsofprod}
    For $i = 1, 2$, let $\mu_i$ denote the Dinh Sibony measure associated to the holomorphic correspondences $F_i$ on $X_{i}$ and let $d_{t_{i}}$ denote the topological degree of $F_i$. Then, we have $\mu_1 \times \mu_2$ to be the corresponding Dinh Sibony measure associated to the product correspondence $F_1 \times F_2$ and hence the support of Dinh Sibony measure of $F_1 \times F_2$ is $\supp (\mu_1) \times \supp (\mu_2)$.
\end{proposition}

\begin{proof}
Let $i = 1, 2$. Since $\mu_i$ is the Dinh Sibony measure associated to the holomorphic correspondences $F_i$, we have, $\lim \limits_{n \, \to \, \infty} \left( \dfrac{1}{d_{t_i}^{n}} \displaystyle \sideset{}{'}  \sum\limits_{x \, \in \, (F_i^n)^{\dagger}(y)} \delta_x\right) = \mu_i$. Now observe that,
\[
  \displaystyle \sideset{}{'}  \sum\limits_{(x_1,x_2) \, \in \, ((F_1 \times F_2)^n)^{{\dagger}}(y_1,y_2)} \frac{\delta_{(x_1,x_2)}}{d_{t_1}^{n}.d_{t_2}^{n}}  = \displaystyle \sideset{}{'}  \sum\limits_{x_1 \, \in \, (F_1^n)^{\dagger}(y_1)} \frac{\delta_{x_1}}{d_{t_1}^{n}}   \ \ \ \times  \displaystyle \sideset{}{'}  \sum\limits_{x_2 \, \in \, (F_2^n)^{\dagger}(y_2)} \frac{\delta_{x_2}}{d_{t_2}^{n}} 
\]
Taking limit on both sides results in,
\[
 \lim \limits_{n \, \to \, \infty} \displaystyle \sideset{}{'}  \sum\limits_{(x_1,x_2) \, \in \, ((F_1 \times F_2)^n)^{{\dagger}}(y_1,y_2)} \frac{\delta_{(x_1,x_2)}}{d_{t_{1}}^{n}.d_{t_{2}}^{n}}  = \left( \lim \limits_{n \, \to \, \infty}  \displaystyle \sideset{}{'}  \sum\limits_{x_1 \, \in \, (F_1^n)^{\dagger}(y_1)} \frac{\delta_{x_1}}{d_{t_1}^{n}}  \right) \times \left(\lim \limits_{n \, \to \, \infty} \displaystyle \sideset{}{'}  \sum\limits_{x_2 \, \in \, (F_2^n)^{\dagger}(y_2)} \frac{\delta_{x_2}}{d_{t_2}^{n}} \right)
 \]
which gives the Dinh Sibony measure of $F_1 \times F_2$ to be $\mu_1 \times \mu_2$. It is now a routine check to observe that the support of Dinh Sibony measure of $F_1 \times F_2$ is $\supp (\mu_1) \times \supp (\mu_2)$.
\end{proof}

The support of the Dinh Sibony measure plays a key role in the study of entropy and pressure related to holomorphic correspondences - much like the analogous role played by the Julia set in case of rational maps on $\widehat{\mathbb{C}}$. A portion of Section \ref{sec:five} deals with the same. Before that, in the next section we introduce the basics of thermodynamic formalism for holomorphic correspondences and prove some results for this theory related to product correspondences.

\section{Thermodynamic formalism for product correspondences}
\label{sec:4}

We first introduce the basic definitions of entropy related to holomorphic correspondences. We follow the ideas mentioned in \cite{Sub2026}. Let $\Gamma$ be a holomorphic correspondence on a compact K\"{a}hler manifold $X$. Let $\nu \in \mathscr{S}^{\Gamma}$ be such that $\nu = (\Pi_{0}^{\Gamma})_{*}\mu$ where $\mu \in \mathscr{M} (\mathscr{P}^{\Gamma} \left( X \right), \sigma^{\Gamma})$ and also consider a finite $\nu$ - measurable partition $\mathcal{Q} = \{Q_1, Q_2, \cdots, Q_s\}$ of $X$. Let $\mathcal{R}_{k}$ and $\mathcal{S}_{k}$ denote the finite $\mu$ - measurable partitions of $\mathscr{P}^{\Gamma} \left( X \right)$ given by $\left\{(\Pi_{k}^{\Gamma})^{-1} (Q_1), (\Pi_{k}^{\Gamma})^{-1} (Q_2), \cdots (\Pi_{k}^{\Gamma})^{-1} (Q_s)\right\}$ and $\left\{(\text{Proj}_{k}^{\Gamma})^{-1} (\{1\}), (\text{Proj}_{k}^{\Gamma})^{-1} (\{2\}), \cdots, (\text{Proj}_{k}^{\Gamma})^{-1} (\{M\})\right\}$ respectively. For $k \ge 0$, consider the partition $\mathcal{Q}_{k}^{\Gamma} := \mathcal{R}_{k} \vee \mathcal{S}_{k+1}$. The entropy of the correspondence $\Gamma$ for the partition $\mathcal{Q}$ and with respect to the pair of measures $(\nu, \mu)$ is defined as, $H_{(\nu, \mu)} (\mathcal{Q}, \Gamma) := \lim \limits_{n \, \to \, \infty} \dfrac{1}{n} H_{\mu} \left(\bigvee_{p = 0}^{n - 1} \mathcal{Q}_{p}^{\Gamma} \right),$
where $H_{\mu} (\mathcal{P})$ denotes the usual notion of entropy of the finite $\mu$ - measurable partition $\mathcal{P}$ with respect to the measure $\mu$. 
We further define the intermediate entropy of $\Gamma$ with respect to the pair of measures $(\nu, \mu)$ as $h_{(\nu, \mu)} (\Gamma) = \sup \limits_{\mathcal{Q}} H_{(\nu, \mu)} (\mathcal{Q}, \Gamma),$ where supremum runs over all the finite $\nu$ measurable partitions of $X$. Finally we are ready to define the measure theoretic entropy of the correspondence $\Gamma$ with respect to the measure $\nu$. 

\begin{definition} 
\label{defn:entropycorr}
Let $\Gamma$ be a holomorphic correspondence defined on $X$ and let $\nu \in \mathscr{S}^{\Gamma}$. Then, the \emph{measure theoretic entropy} of the holomorphic correspondence $\Gamma$ with respect to the measure $\nu$ is defined as $h_{\nu} \left( \Gamma \right) = \sup \left\{ h_{\left( \nu,\, \mu \right)} \left( \Gamma \right) : \mu \in \mathscr{M} \left( \mathscr{P}^{\Gamma} \left( X \right), \sigma^{\Gamma} \right)\ \text{such that} \ \nu = \left( \Pi_{0}^{\Gamma} \right)_{*} \mu \right\}.$ 
\end{definition} 

The following is a theorem that explicitly finds the value of $h_{(\nu, \mu)}$ for appropriately chosen measures. This result considerably simplifies the calculations related to $h_{\nu} (\Gamma)$. 

\begin{theorem}\cite{Sub2026}
    \label{thm:entropeq}
    Let $\Gamma$ be a holomorphic correspondence defined on $X$ and $\sigma^{\Gamma}$ denote the shift map on the orbit space $\mathscr{P}^{\Gamma} \left( X \right)$. Then for all pair of measures $(\nu, \mu)$ with $\nu \in \mathscr{S}^{\Gamma}$ and $\mu \in \mathcal{M} (\mathscr{P}^{\Gamma} \left( X \right), \sigma^{\Gamma})$ such that $(\Pi_{0}^{\Gamma})_{*} \mu = \nu$, we have $h_{(\nu, \mu)} (\Gamma) = h_{\mu} (\sigma^{\Gamma})$, where $h_{\mu} (\sigma^{\Gamma})$ denotes the measure theoretic entropy of the map $\sigma^{\Gamma}$ with respect to the measure $\mu$. 
\end{theorem}

We now state a theorem from \cite{walters_book} that comes in handy while proving  Theorem \ref{thm:entropycorrineq}. 

\begin{theorem}
\cite{walters_book}
    \label{thm:entprodmaps}
    Fix $i = 1, 2$. Let $T_i$ be continuous maps on the compact metric space $X_i$ and let $\mu_i$ be a $T_i$ invariant Borel probability measure. Then we have, $h_{\mu_1 \times \mu_2} ((T_1 , T_2)) = h_{\mu_1} (T_1) + h_{\mu_2} (T_2),$ where $h_{\mu_i} (T_i)$ denotes the measure theoretic entropy of $T_i$ with respect to $\mu_i$ and $h_{\mu_1 \times \mu_2} ((T_1 , T_2))$ denotes the measure theoretic entropy of $(T_1 , T_2)$ with respect to $\mu_1 \times \mu_2$.
\end{theorem}

We are now ready to state and prove the theorem on the relation between the measure theoretic entropies of the product correspondence and the individual correspondences. 

\begin{theorem}
\label{thm:entropycorrineq}
    For $i = 1, 2$, let $\Gamma_{i}$ be holomorphic correspondences on $X_i$. For a measure $\nu_i \in \mathscr{S}^{\Gamma_{i}}$, let $h_{\nu_i} \left( \Gamma_i \right)$ denote the measure theoretic entropy of the holomorphic correspondences $\Gamma_i$ with respect to the measure $\nu_i$. Then, we have
    \[h_{\nu_1}\left( \Gamma_1 \right) + h_{\nu_2}\left( \Gamma_2 \right) \leq h_{(\nu_1 \times \nu_2)}\left( \Gamma_1 \times \Gamma_2\right).\]
\end{theorem}

\begin{proof}
Let $\mu_1 \in \mathscr{M}\left(\mathscr{P}^{\Gamma_1}{(X_1)}, \sigma^{\Gamma_{1}} \right)$ and $\mu_2 \in \mathscr{M}\left(\mathscr{P}^{\Gamma_2}{(X_2)}, \sigma^{\Gamma_2} \right)$ satisfy $(\Pi_{0}^{\Gamma_{1}})_{*} (\mu_1) = \nu_1$ and $(\Pi_{0}^{\Gamma_{2}})_{*} (\mu_2) = \nu_2$. Proposition \ref{prop:SGamma} and Theorem \ref{thm:entropeq} together give $h_{\left(\nu_1 \times \nu_2, \, \Phi_*(\mu_1 \times \mu_2)\right)}(\Gamma_{1} \times \Gamma_{2}) = h_{\Phi_*(\mu_1 \times \mu_2)}\left(\sigma^{\Gamma_1 \times \Gamma_2}\right).$ \\

Now using the Definition \ref{defn:entropycorr}, Lemma \ref{lem:MPS}, Remark \ref{rem:entpreseq} and Theorem \ref{thm:entprodmaps} we have,
\begin{align*}
h_{(\nu_1 \times \nu_2)}(\Gamma_1 \times \Gamma_2) 
&= \sup \Big\{ h_{(\nu_1 \times \nu_2, \, \mu)}(\Gamma_1 \times \Gamma_2) : \mu \in \mathscr{M}\left(\mathscr{P}^{{\Gamma_1}\times{\Gamma_2}}(X_1 \times X_2), \sigma^{\Gamma_1 \times \Gamma_2}\right), \Big. \\
& \Big. \hspace{0.4cm}\text{ with } (\Pi_0^{\Gamma_1 \times \Gamma_2})_{*} \mu = \nu_1 \times \nu_2 \Big\} \\
&= \sup \Big\{ h_{(\nu_1 \times \nu_2, \, \Phi_* \mu')}(\Gamma_1 \times \Gamma_2) : \mu' \in \mathscr{M}\left(\mathscr{P}^{\Gamma_1}{(X_1)} \times \mathscr{P}^{\Gamma_2}{(X_2)}, (\sigma^{\Gamma_1}, \sigma^{\Gamma_2})\right), \Big. \\
& \Big. \hspace{0.4cm}\text{ with } (\Pi_0^{\Gamma_1 \times \Gamma_2})_{*} (\Phi_* \mu') = \nu_1 \times \nu_2 \Big\} \\
& \ge \sup \Big\{ h_{(\nu_1 \times \nu_2, \, \Phi_*(\mu_1 \times \mu_2))}(\Gamma_1 \times \Gamma_2) : \mu_1 \in \mathscr{M}\left(\mathscr{P}^{\Gamma_1}{(X_1)}, \sigma^{\Gamma_1}\right), \Big. \\
& \Big. \hspace{0.4cm}\mu_2 \in \mathscr{M}\left(\mathscr{P}^{\Gamma_2}{(X_2)}, \sigma^{\Gamma_2}\right), \text {with } (\Pi_0^{\Gamma_1})_{*} (\mu_1) = \nu_1 \text{ and } (\Pi_0^{\Gamma_2})_{*} \mu_2 = \nu_2 \Big\} \\
&= \sup \Big\{ h_{(\nu_1, \mu_1)}(\Gamma_1) + h_{(\nu_2, \mu_2)}(\Gamma_2) : \mu_1 \in \mathscr{M}\left(\mathscr{P}^{\Gamma_1}{(X_1)}, \sigma^{\Gamma_1}\right), \Big.\\
& \Big. \hspace{0.4cm}\mu_2 \in \mathscr{M}\left(\mathscr{P}^{\Gamma_2}{(X_2)}, \sigma^{\Gamma_2}\right), \text {with } (\Pi_0^{\Gamma_1})_{*} (\mu_1) = \nu_1 \text{ and } (\Pi_0^{\Gamma_2})_{*} \mu_2 = \nu_2 \Big\} \\
&= h_{\nu_1}(\Gamma_1) + h_{\nu_2}(\Gamma_2)
\end{align*}
\end{proof}

We now proceed to studying the concept of pressure in the case of holomorphic correspondences. Pressure is defined using some special sets called separated and spanning sets, which we shall now introduce.

\begin{definition}
Let $\epsilon > 0$ and $n \in \mathbb{N}$.
\begin{enumerate}
    \item A subset $\mathcal{F}$ of $\mathscr{P}_{n}^{\Gamma} (X)$ is said to be $(n, \epsilon)$ separated if for any distinct points $\mathfrak{X}_{n}^{+} (x_{0}; \boldsymbol{\alpha})_{\boldsymbol{j}}$ and $\mathfrak{X}_{n}^{+} (y_{0}; \boldsymbol{\beta})_{\boldsymbol{k}}$ in $\mathcal{F}$ , we have
     \begin{eqnarray*}
        \text{either} \ \  & d_{X} \left(\Pi_{(r,n)}^{\Gamma}  \left( \mathfrak{X}^{+}_{n} \left( x_{0}; \boldsymbol{\alpha} \right)_{\boldsymbol{j}} \right), \Pi_{(r,n)}^{\Gamma}  \left( \mathfrak{X}^{+}_{n} \left( y_{0}; \boldsymbol{\beta} \right)_{\boldsymbol{k}} \right)  \right)  > \epsilon, \ \text{for some} \ 0 \le r \le n  \\
        \text{or} \ \ & \text{Proj}_{(j,n)}^{\Gamma} \left(\mathfrak{X}^{+}_{n} \left( x_{0}; \boldsymbol{\alpha} \right)_{\boldsymbol{j}}\right) \neq  \text{Proj}_{(j,n)}^{\Gamma} \left(\mathfrak{X}^{+}_{n} \left( y_{0}; \boldsymbol{\beta} \right)_{\boldsymbol{k}}\right), \ \text{for some} \  1 \le j \le n.
    \end{eqnarray*}
    \item A subset $\mathcal{G}$ of $\mathscr{P}_{n}^{\Gamma} (X)$ is said to be $(n, \epsilon)$ spanning set of $\mathscr{P}_{n}^{\Gamma} (X)$ if for every $\mathfrak{X}^{+}_{n} \left( x_{0}; \boldsymbol{\alpha} \right)_{\boldsymbol{j}} \in \mathscr{P}_{n}^{\Gamma} (X)$, there exists $y_0 \in X$ such that $\mathfrak{X}^{+}_{n} \left( y_{0}; \boldsymbol{\alpha} \right)_{\boldsymbol{k}} \in \mathcal{G}$ satisfying 
    \[
    d_{X} \left(\Pi_{(r,n)}^{\Gamma}  \left( \mathfrak{X}^{+}_{n} \left( x_{0}; \boldsymbol{\alpha} \right)_{\boldsymbol{j}} \right), \Pi_{(r,n)}^{\Gamma}  \left( \mathfrak{X}^{+}_{n} \left( y_{0}; \boldsymbol{\alpha} \right)_{\boldsymbol{k}} \right)  \right)  < \epsilon, \ \text{for all} \ 0 \le r \le n.
    \]
    \end{enumerate}
\end{definition}

Before defining the pressure of a continuous map with respect to the given holomorphic correspondence, we first give a few notations. For $f \in \mathcal{C} (X, \mathbb{R})$ and $n \in \mathbb{N}$, we define,
\[
\mathcal{Z}_{n} (f) (\mathfrak{X}_{n}^{+} (x_{0}; \boldsymbol{\alpha})_{\boldsymbol{j}}) = \sum \limits_{j \, = \, 0}^{n - 1} f \left(\Pi_{(r,n)}^{\Gamma}  \left( \mathfrak{X}^{+}_{n} \left( x_{0}; \boldsymbol{\alpha} \right)_{\boldsymbol{j}} \right) \right),
\]
for any $\mathfrak{X}_{n}^{+} (x_{0}; \boldsymbol{\alpha})_{\boldsymbol{j}} \in \mathscr{P}_{n}^{\Gamma} (X)$. Now for $\epsilon > 0$, we define
\begin{eqnarray*}
    \mathscr{Q}_{n}^{\Gamma} (f, \epsilon) & = & \inf \limits_{\mathcal{G}} \left \{ \sum \limits_{\mathfrak{X}_{n}^{+} (x_{0}; \boldsymbol{\alpha})_{\boldsymbol{j}} \in \mathcal{G}} e^{\mathcal{Z}_{n} (f) (\mathfrak{X}_{n}^{+} (x_{0}; \boldsymbol{\alpha})_{\boldsymbol{j}})} \right \}, \\
    & & \text{where the infimum is over all possible} \ (n, \epsilon) \ \text{spanning sets of} \ \mathscr{P}_{n}^{\Gamma} (X). \\ 
    \mathscr{R}_{n}^{\Gamma} (f, \epsilon) & = & \sup \limits_{\mathcal{F}} \left \{ \sum \limits_{\mathfrak{X}_{n}^{+} (x_{0}; \boldsymbol{\alpha})_{\boldsymbol{j}} \in \mathcal{F}} e^{\mathcal{Z}_{n} (f) (\mathfrak{X}_{n}^{+} (x_{0}; \boldsymbol{\alpha})_{\boldsymbol{j}})} \right \}, \\  
    & & \text{where the supremum is over all possible} \ (n, \epsilon) \ \text{separated sets of} \ \mathscr{P}_{n}^{\Gamma} (X). 
\end{eqnarray*}

\begin{definition}
\label{def:pressure corr}
    Let $\Gamma$ be a holomorphic correspondence on a K\"{a}hler manifold $X$ and let $f \in \mathcal{C} (X, \mathbb{R})$. We define the pressure of $f$ with respect to $\Gamma$ to be 
    \[
    {\rm Pr} (\Gamma, f) = \lim_{\epsilon \, \to \, 0} \left ( \limsup \limits_{n \, \to \, \infty} \frac{1}{n} \log \mathscr{Q}_{n}^{\Gamma} (f, \epsilon) \right) = \lim_{\epsilon \, \to \, 0} \left ( \limsup \limits_{n \, \to \, \infty} \frac{1}{n} \log \mathscr{R}_{n}^{\Gamma} (f, \epsilon) \right).
    \]
\end{definition}

Our next aim is to prove a result on topological pressure of product correspondence analogous to Theorem \ref{thm:entropycorrineq}. The following is a theorem from \cite{Subith_ruelle} which comes handy for that task. Note that in \cite{Subith_ruelle} the authors work with the domain being $\widehat{\mathbb{C}}$. However, the techniques of the proof of the statement remains the same even if we take the domain to be a compact K\"{a}hler manifold $X$. 

\begin{theorem}
\label{thm:pressureGamma}
Let $\Gamma$ be a holomorphic correspondence on $X$. For every continuous function $f \in \mathcal{C} \left( X, \mathbb{R} \right)$, we have ${\rm Pr} (\Gamma, f) = {\rm Pr} (\sigma^{\Gamma}, f \circ \Pi_{0}^{\Gamma})$, where  ${\rm Pr} (\sigma^{\Gamma}, f \circ \Pi_{0}^{\Gamma})$ is the topological pressure of $f \circ \Pi_{0}^{\Gamma} \in \mathcal{C} (\mathscr{P}^{\Gamma} \left( X \right), \mathbb{R})$, with respect to the shift map $\sigma^{\Gamma}$. 
\end{theorem}

It turns out that, using the above theorem, the relation between pressure associated to the product correspondence $\Gamma_1 \times \Gamma_2$ with respect to an appropriate potential has a nice characterisation with respect to the pressure of the correspondences $\Gamma_i$ with respect to the potential $f_i$, where $i = 1,2$. In the upcoming theorem, we use the following notation. For functions $f_i \in \mathcal{C} (X_i, \mathbb{R})$, we denote by $f_1 \times f_2 \in \mathcal{C} (X_1 \times X_2, \mathbb{R})$, the potential defined as $(f_1 \times f_2) (x_1, x_2) = f_1(x_1) + f_2(x_2)$. We know that for maps $T_i$ defined on $X_i$, we the equality related to pressures given by ${\rm Pr} ((T_1, T_2), f_1 \times f_2) = {\rm Pr} (T_1, f_1) + {\rm Pr} (T_2, f_2)$. We use this result to prove the following theorem.

\begin{theorem}
\label{thm:presseqcorr}
    For $i = 1, 2$, let $\Gamma_{i}$ be holomorphic correspondences on $X_i$ and let $f_i \in \mathcal{C} (X_i, \mathbb{R})$. Then we have, ${\rm Pr} (\Gamma_1 \times \Gamma_2, f_1 \times f_2)={\rm Pr} (\Gamma_1, f_1)+{\rm Pr} (\Gamma_2, f_2)$.
\end{theorem}

\begin{proof}
    Using Theorem \ref{thm:pressureGamma} and Remark \ref{rem:entpreseq}. We have,
\begin{eqnarray*}
    {\rm Pr} (\Gamma_1 \times \Gamma_2, f_1 \times f_2)&=&{\rm Pr}\left(\sigma^{\Gamma_1 \times \Gamma_2}, (f_1 \times f_2)\circ \Pi_0^{\Gamma_1 \times \Gamma_2}\right)\\
    &=& {\rm Pr}\left( \left(\sigma^{\Gamma_1}, \sigma^{\Gamma_2}\right), (f_1 \times f_2)\circ \Pi_0^{\Gamma_1 \times \Gamma_2} \circ \Phi \right)\\
    &=& {\rm Pr}\left( \left(\sigma^{\Gamma_1}, \sigma^{\Gamma_2}\right), (f_1 \circ \Pi_0^{\Gamma_1}) \times (f_2 \circ \Pi_0^{\Gamma_2}) \right)\\
    &=& {\rm Pr} (\Gamma_1, f_1)+{\rm Pr} (\Gamma_2, f_2),
\end{eqnarray*}
which proves the theorem.
\end{proof}

\section{Some families of correspondences}
\label{sec:five}

In this section we explore various families of holomorphic correspondences possessing some interesting properties and prove results on whether the product correspondence preserves these properties provided the individual correspondences comprising the product have this property. The underlying idea that is studied in detail in this section is related to various notions of expanding criteria. We study each notion and also make comments on how they are different from the other. We also give some remarks about how these concepts work in the case of maps and what changes are observed in the case of holomorphic correspondences. An important component of this section is the result regarding uniqueness of the equilibrium state, that serves as a nice application of expanding behaviour in correspondences with regards to thermodynamic formalism. First we give a brief remark on restrictive holomorphic correspondences.

\subsection{Restrictive correspondences}

The concept of what we define below as restrictive correspondences has its origins in the work \cite{quad_corr} by Bullett and was further studied in detail in \cite{rest_corr}. The authors in that paper defined this concept only for polynomial correspondences. Here we provide a natural extension of that study to the general holomorphic correspondence setting. Let $F$ be a holomorphic correspondence on a compact K\"{a}hler manifold $X$ with topological degree $d_t$ and forward degree $d_f$. For a generic point $x \in X$, let $F(x) = \{ y_1 (x), y_2 (x), \cdots, y_{d_{f}}(x) \}$.

\begin{definition}
\label{defn:restcorr}
    A holomorphic correspondence $F$ on a compact K\"{a}hler manifold $X$ is said to be \emph{restrictive} if for all $x \in X$, we have $F^{\dagger} \left(y_{j} (x)\right)$ to be the same set for each $1 \le j \le d_f$. 
\end{definition}
A simple calculation shows that whenever $F$ is restrictive we have $F^{\dagger}$ to be restrictive as well. Suppose $F_1$ and $F_2$ are restrictive holomorphic correspondences on $X_1$ and $X_2$ with forward degrees $d_{f_{1}}$ and $d_{f_{2}}$ respectively and with topological degrees to be $d_{t_{1}}$ and $d_{t_{2}}$ respectively. Then, it is easy to see that $F_1 \times F_2$ is also restrictive. Indeed if $(x_1, x_2) \in X_1 \times X_2$ and if we have $F_1(x_1) = \{y_{1} (x_1), y_{2} (x_1), \cdots, y_{d_{t_{1}}} (x_1) \}$ and $F_2(x_2) = \{y_{1} (x_2), y_{2} (x_2), \cdots, y_{d_{t_{2}}} (x_2) \}$, then for $1 \le j \le d_{f_{1}}$ and $1 \le k \le d_{f_{2}}$, we have by definition, $(F_1 \times F_2)^{\dagger} ((y_{j} (x_1), y_{k} (x_2))) = F_1^{\dagger} (y_{j} (x_1)) \times F_2^{\dagger} (y_{k} (x_2)).$ The result now follows as we have $F_1^{\dagger} (y_{j} (x_1))$ to be the same set for each $1 \le j \le d_{f_{1}}$ and $F_2^{\dagger} (y_{k} (x_2))$ is the same set for each $1 \le k \le d_{f_{2}}$. 

\subsection{Distance expanding correspondences and unique equilibrium state}

We now turn our attention to the main concept of this section namely, studying distance expanding holomorphic correspondences.

\begin{definition}
\label{defn:distexp}
    A holomorphic correspondence defined on a compact K\"{a}hler manifold $X$ is said to be \emph{distance expanding} if $\exists \, \lambda > 1$ and $\delta > 0$ such that for every distinct $x, y \in X$ such that $d_{X} (x, y) < \delta$ we have $\inf \left\{ d_{X} (x', y') : x' \in F(x), \, y' \in F(y) \right\} > \lambda \, d_{X} (x, y)$.
\end{definition}

\begin{definition}
\label{defn:fwdexps}
    A holomorphic correspondence defined on a compact K\"{a}hler manifold $X$ is said to be \emph{forward expansive} if $\exists \, \delta > 0$ such that for all distinct $x, y \in X$, $\exists \, n \geq 0$ satisfying $\inf \{d_{X}(x',y') \, |\, x' \in F^{n}(x) \, \, \text{and} \, \, y' \in F^{n}(y)\} > \delta$.
\end{definition}

\begin{proposition}
\label{prop:proddistexp}
    Let $F_1$ and $F_2$ be holomorphic correspondences on compact K\"{a}hler manifolds $X_1$ and $X_2$ respectively.
    \begin{enumerate}
        \item Suppose $F_1$ and $F_2$ are distance expanding, then so is $F_1 \times F_2$.
        \item Suppose $F_1$ and $F_2$ are forward expansive, then so is $F_1 \times F_2$.
    \end{enumerate}
\end{proposition}

\begin{proof} Let $F_{j}$ be a holomorphic correspondence defined on $X_{j}$, where $j = 1,2$.
\begin{enumerate}
    \item  Consider the constants $\lambda_j$ and $\delta_{j}$ associated to the correspondence $F_j$, guaranteed by the hypothesis and satisfying the condition in Definition \ref{defn:distexp}. Set $\lambda = \min \{\lambda_1, \lambda_2\}$, $\delta = \min\{\delta_1, \delta_2\}$ and consider points with $d_{X_1 \times X_2} ((x_1, x_2), (y_1, y_2)) < \delta$.
    
    Then, by hypothesis, for such points, we have, $d_{X_1} (x_1', y_1') \geq \lambda d_{X_1} (x_1, y_1)$ for all $x_1' \in F_{1} (x_1)$ and $y_1' \in F_{1} (y_1)$. Similarly, we have $d_{X_2} (x_2', y_2') \geq \lambda d_{X_2} (x_2, y_2)$ for all $x_2' \in F_{2} (x_2)$ and $y_2' \in F_{2} (y_2)$. Hence for all $(x_1', x_2') \in (F_1 \times F_2) (x_1, x_2)$ and $(y_1', y_2') \in (F_1 \times F_2) (y_1, y_2)$ we have, $\max \{d_{X_1} (x_1', y_1'), d_{X_{2}}(x_2', y_2')\} \geq \lambda \max \{d_{X_1} (x_1, y_1), d_{X_{2}}(x_2, y_2)\}$
    and consequently, 
    \begin{eqnarray*}
        & & \inf \left \{ d_{X_1 \times X_2} ((x_1', x_2'), (y_1', y_2')) : (x_1', x_2') \in (F_1 \times F_2) (x_1, x_2), \, (y_1', y_2') \in (F_1 \times F_2) (y_1, y_2) \right\} \\
        &  & \hspace{+50mm} \geq  \lambda d_{X_1 \times X_2} ((x_1, x_2), (y_1, y_2)).
    \end{eqnarray*}

    \item Using forward expansiveness of $F_{i}$, we get $\delta_{i} > 0$ such that the condition in the Definition \ref{defn:fwdexps} is satisfied. Set $\delta := \min \{\delta_{1}, \delta_{2}\}$. 
    
    Consider $(x_{1}, x_{2}), (y_{1}, y_{2}) \in X_{1} \times X_{2}$ with $(x_{1}, x_{2}) \neq (y_{1}, y_{2}) $. Suppose $x_{1} \neq y_{1}$. Then, there exists $n \geq 0$ such that, $\inf \{d_{X_{1}}(x_{1}',y_{1}') \, |\, x_{1}' \in F_{1}^{n}(x_{1}) \, \, \text{and} \, \, y_{1}' \in F_{1}^{n}(y_{1})\} > \delta.$ Hence, for any $(x_{1}', x_{2}') \in (F_1 \times F_2)^{n}(x_{1}, x_{2})$, and for $(y_{1}', y_{2}') \in (F_1 \times F_2)^{n}(y_{1}, y_{2})$,
    \[
    d_{X_{1}\times X_{2}}\left((x_{1}', x_{2}'),(y_{1}', y_{2}')\right) =  \max \{ d_{X_{1}}(x_{1}', y_{1}'), d_{X_{2}}(x_{2}', y_{2}')\} \geq  d_{X_{1}}(x_{1}', y_{1}') >  \delta 
    \]
    This implies
    \[
    \inf \{d_{X_{1}\times X_{2}}\left((x_{1}', x_{2}'),(y_{1}', y_{2}')\right) \, | \, (x_{1}', x_{2}') \in (F_1 \times F_2)^{n}(x_{1}, x_{2}), (y_{1}', y_{2}') \in (F_1 \times F_2)^{n}(y_{1}, y_{2})\} > \delta.
    \]
    For points $(x_{1}, x_{2}), (y_{1}, y_{2}) \in X_{1} \times X_{2}$ with $x_2 \ne y_2$, one uses the fact that $F_2$ is forward expansive to get the above inequality. 
    \end{enumerate}
\end{proof}

An important application of distance expanding maps with regards to thermodynamic formalism is about the uniqueness of equilibrium state associated to a H\"{o}lder continuous potential. Let $F$ be a holomorphic correspondence on $X$ and let $\Omega_{{\rm DS}}$ denote the support of the Dinh Sibony measure of $F$. In the remainder of this subsection, we study the dynamics of holomorphic correspondences restricted to their support of Dinh Sibony measures. The following variational principle from \cite{Sub2026} was the first step towards studying the concept of unique equilibrium state. While they proved it for correspondences defined on $\widehat{\mathbb{C}}$, it is easy to see that this result holds in the general setting as well.

\begin{theorem}
\label{thm:varprincicorr}
    Let $F$ be a holomorphic correspondence on $X$ and let $f \in \mathcal{C} (\Omega_{{\rm DS}}, \mathbb{R})$. Then, we have ${\rm Pr} (\Gamma, f) = \sup \limits_{\nu \in \mathscr{S}^{\Gamma}} \left \{ h_{\nu} (\Gamma) + \displaystyle\int f \, \mathrm{d}\nu  \right \}.$
\end{theorem} 

Any measure $\nu \in \mathscr{S}^{\Gamma}$ where the above supremum is attained is called an equilibrium measure (or equilibrium state) for $f$. Following the work \cite{Sub2026}, the results proved in \cite{uniqueequibcorr} show that under the dynamics of a holomorphic correspondence $F$ which is distance expanding on $\Omega_{{\rm DS}}$, any H\"{o}lder continuous potential on $\Omega_{{\rm DS}}$ has a unique equilibrium state. We now investigate this problem in case of product of two holomorphic correspondences. Before that, we formally state the main result from \cite{uniqueequibcorr}. 

\begin{theorem}
\label{thm:unqiequibcorr}
    Let $F$ be a holomorphic correspondence which is distance expanding on $\Omega_{{\rm DS}}$. Then any H\"{o}lder continuous function $f \in \mathcal{C} (\Omega_{{\rm DS}}, \mathbb{R})$ has a unique equilibrium state.
\end{theorem}

Suppose $F_1$ and $F_2$ are holomorphic correspondences on $X_1$ and $X_2$ for which the support of the Dinh Sibony measure exists and are denoted by $\Omega_{{\rm DS}}^{1}$ and $\Omega_{{\rm DS}}^{2}$ respectively. Assume $F_i$ is distance expanding on $\Omega_{{\rm DS}}^{i}$. By Proposition \ref{prop:dsofprod}, we have the support of Dinh Sibony measure of $F_1 \times F_2$ to be $\Omega_{{\rm DS}}^{1} \times \Omega_{{\rm DS}}^{2}$. We now restrict the correspondences $F_1, F_2$ and $F_1 \times F_2$ to their respective support of Dinh Sibony measures. The spaces of orbits, measures, etc. are also considered with this set being the domain of the correspondence. \vspace{2mm}

Let $f_1: \Omega_{{\rm DS}}^{1} \to \mathbb{R}$ and $f_2: \Omega_{{\rm DS}}^{2} \to \mathbb{R}$ be H\"{o}lder continuous functions and let $\nu_1 \in \mathscr{S}^{\Gamma_{1}}$ and $\nu_{2} \in \mathscr{S}^{\Gamma_{2}}$ be their respective unique equilibrium states as guaranteed by Theorem \ref{thm:unqiequibcorr}. Consider the function $f_1 \times f_2: \Omega_{{\rm DS}}^{1} \times \Omega_{{\rm DS}}^{2} \to \mathbb{R}$ given by $(f_1 \times f_2) (x_1, x_2) = f_1 (x_1) + f_2 (x_2)$. Then we have the following theorem.

\begin{theorem}
\label{thm:uniqequi}
    The unique equilibrium state associated to the potential $f_1 \times f_2$ under the dynamics of $F_1 \times F_2$ is $\nu_1 \times \nu_2$.
\end{theorem}

\begin{proof}
From Theorem \ref{thm:entropycorrineq}, we have
\[
h_{\nu_1} (\Gamma_{1}) + h_{\nu_{2}} (\Gamma_{2}) + \int \limits_{\Omega_{{\rm DS}}^{1}} f_1 \, \mathrm{d}\nu_1 + \int \limits_{\Omega_{{\rm DS}}^{2}} f_2 \, \mathrm{d}\nu_2 \ \le \ h_{\nu_1 \times \nu_2} (\Gamma_1 \times \Gamma_2) + \int \limits_{\Omega_{{\rm DS}}^{1} \times \Omega_{{\rm DS}}^{2}} f_1 \times f_2 \, \mathrm{d}(\nu_1 \times \nu_2).
\]
Using Theorem \ref{thm:presseqcorr} and the fact that $\nu_1$ and $\nu_2$ are appropriate equilibrium states, we get that ${\rm Pr} (\Gamma_1 \times \Gamma_2, f_1 \times f_2) \le \ h_{\nu_1 \times \nu_2} (\Gamma_1 \times \Gamma_2) + \displaystyle \int \limits_{\Omega_{{\rm DS}}^{1} \times \Omega_{{\rm DS}}^{2}} f_1 \times f_2 \, \mathrm{d}(\nu_1 \times \nu_2).$ Now, an application of Theorems \ref{thm:varprincicorr} and \ref{thm:unqiequibcorr} completes the proof.
\end{proof}

\subsection{Some remarks} 

This subsection focuses on comparing expanding behaviour in case of maps and that of correspondences. We also make a comment on the expansive correspondences condition introduced in \cite{Subith_ruelle} in order to study a Ruelle operator for holomorphic correspondences. \vspace{2mm}

A holomorphic correspondence defined on a compact K\"{a}hler manifold $X$ is said to be \emph{globally distance expanding} if $\exists \, \lambda > 1$  such that for every $x, y \in X$ such that $x \neq y$, we have $\inf \left\{ d_{X} (x', y') : x' \in F(x), \, y' \in F(y) \right\} > \lambda \, d_{X} (x, y)$. It is easy to observe that for a globally distance expanding correspondence, $F(x) \cap F(y) = \emptyset$ for all $x, y \in X$ such that $x \neq y$. This forces $F(x) \cap F(y) = \emptyset$. For $x' \in F(x), \, y' \in F(y)$, applying the global distance expanding property of $F$ gives, $\inf \left\{ d_{X} (x'', y'') : x'' \in F(x'), \, y'' \in F(y') \right\} > \lambda \, d_{X} (x', y') > \lambda^{2} \, d_{X} (x, y).$ Similarly, we can show that $F^n$ is globally distance expanding with expanding constant $\lambda^n$ for any $n \in \mathbb{N}$. Fix any $\delta > 0$ and let $x \ne y$ and choose $n \in \mathbb{N}$ such that, $\lambda^n d(x, y) > \delta$. Thus $\inf \left\{ d_{X} (z', w') : z' \in F^{n}(z), \, w' \in F^{n}(w') \right\} > \lambda^{n} \, d_{X} (z, w) > \delta$. The discussion in this paragraph is summarised as the following result.

\begin{result}
\label{thm:globexp}
    If $F$ is a globally distance expanding holomorphic correspondence on $X$, then so is $F^n$, for any $n \in \mathbb{N}$. Moreover, $F$ is forward expansive.
\end{result}

Typically, in the case of maps, distance expanding (which is defined locally) implies forward expansiveness. But in correspondences, it becomes difficult to do it with the local definition as a point can have multiple forward images. Thus, finding local neighbourhoods around each of the forward images at each stage becomes difficult. This is one of the motivations for us to consider globally distance expanding correspondences. Further, note that in the case of the domain $X$ being the Riemann sphere $\widehat{\mathbb{C}}$, the correspondences which are distance expanding with $\delta > 2$, where $\delta$ denotes the radius of local neighbourhood as mentioned in Definition \ref{defn:distexp}, automatically becomes globally distance expanding on $\widehat{\mathbb{C}}$ as the spherical metric therein is bounded by $2$. Also, as observed before, if $F$ is globally distance expanding, then for any $x \neq y$, we have $F(x) \cap F(y) = \emptyset$. This is the analogue of the fact that if $T$ is a distance expanding map, then it is locally injective. \vspace{2mm}

Another notion of measuring how fast the orbits move under the dynamics of a correspondence was studied in \cite{Subith_ruelle}. The authors of that work studied what they called as expansive correspondences. One of the main objectives of the authors in \cite{Subith_ruelle} while studying expansive correspondences was to explore a Ruelle operator in the setting of holomorphic correspondences. A holomorphic correspondence $\Gamma$ on $X$ with a representation as given in Equation \eqref{eqn:holocorr} is said to be an \emph{expansive correspondence} on $X$ if there exists $\lambda > 1$ such that for any distinct $x, y \in X$ and for any $1 \le j \le N$ with $(x_{-1}, x) \in \Gamma_{j}$, there exists $y_{-1} \in X$ with $(y_{-1}, y) \in \Gamma_{j}$ such that $\lambda \, d_{X} (x_{-1}, y_{-1}) < d_{X} (x, y)$. \vspace{2mm}

Suppose $\Gamma$ is a globally distance expanding correspondence, then for any $x, y \in X$, any fixed $1 \le j \le N$ and any point $(x_{-1}, x) \in \Gamma_{j}$, if we consider $y_{-1} \in X$ to be such that $(y_{-1}, y) \in \Gamma_{j}$, using the definition, there exists $\lambda > 1$ such that,
\[
    d_{X} (x, y) \ge \inf \{ d_{X} (x', y'): x' \in F(x_{-1}), \, y \in F(y_{-1}) \} > \lambda \, d(x_{-1}, y_{-1}),
\]
which gives a simple proof of the following result.

\begin{result}
\label{rslt:GDE-Exp}
    If $F$ is a globally distance expanding correspondence on $X$, then it is an expansive correspondence.
\end{result}

One possible application of the above result is that we can apply the theory of Ruelle operators developed in \cite{Subith_ruelle} to the globally distance expanding correspondences too.

\nocite{*}
\printbibliography[heading=bibintoc,title=References]
\thispagestyle{empty}

\end{document}